\documentclass[letterpaper,11pt]{amsart}
\usepackage{amsmath}
\usepackage{amssymb}
\usepackage{amsfonts}
\usepackage{amsthm}
\def\inv{^{-1}}
 \DeclareMathOperator{\xc}{XC}
\DeclareMathOperator{\Ext}{Ext} \DeclareMathOperator{\Hom}{Hom}

\newcommand{\bsl}{\begin{slide}{}}

\newcommand{\h}{\mathcal{H}}

\newcommand{\A}{\mathcal{A}}

\newcommand{\X}{\mathcal{X}}

\newcommand{\Cal}[1]{\mathcal #1}

\def\sbr #1.{^{[#1]}}
\def\sfl #1.{^{\lfloor #1\rfloor}}

\newcommand\jsh{\rm J^\sharp}
\newcommand\gb{\mathfrak g_\bullet}
\newcommand\half{\frac{1}{2}}

\def\what{\widehat}
\def\inv{^{-1}}
\def\?{{\bf{??}}}

\def\HH{\mathbb H}

\def\A{\Bbb A}
\def\X{\Cal X}

\def\C{\Bbb C}

\def\Spec{\text{\rm Spec} }
\def\ad{\text{\rm ad}}

\def\O{\Cal O}

\def\Sym{\textrm{Sym}}
\def\id{\text{id}}

\def\g{\frak g}

\def\gl{\frak{gl}}

\def\h{\frak h}
\def\rs{{\rho\sigma}}
\def\ab{{\alpha\beta}}
\def\bg{{\beta\gamma}}
\def\ga{{\gamma\alpha}}
\def\abg{{\alpha\beta\gamma}}
\def\m{\frak m}

\def\hom{\frak {hom}}

\def\1/2{\frac{1}{2}}

\def\U{\frak U}

\def\2{{[2]}}
\def\l{\ell}
\def\nl{\newline}
\def\n{\mathcal{N}}
\def\t{\mathcal{T}}
\def\hom{\mathcal{H}\mathit{om}}

\def\<{\langle}
\def\>{\rangle}
\def\sela{SELA\ }
\def\endo{\mathcal{E}\mathit{nd}}

\def\beq{\begin{equation}}
\def\eeq{\end{equation}}
\newcommand{\mex}[1]{\begin{example}\begin{rm}#1
\end{rm}
\end{example}}
\def\eex{\end{rm}\end{example}}

\newtheorem{thm}{Theorem}[section]

\newtheorem{lem}[thm]{Lemma}

\newtheorem{rem}[thm]{Remark}
\newtheorem{example}[thm]{Example}
\begin{document}
\title{ Bernoulli numbers\\ and deformations of schemes and maps}

\author{Ziv Ran}
\thanks{ Research supported in part by NSA grant
 H98230-05-1-0063; v.060704}
\address{University of California,
 Riverside}
\email{ziv.ran@ucr.edu}
\date{\today}\begin{abstract} We introduce a notion of Jacobi-Bernoulli
cohomology associated to a semi-simplicial Lie algebra (SELA). For
an algebraic scheme $X$ over $\C$, we construct a tangent SELA
$\t_X$ and show that the Jacobi-Bernoulli cohomology of $\t_X$ is
related to infinitesimal deformations of $X$.
\end{abstract}

\maketitle \setcounter{section}{-1} \section{Overview} The 'usual'
deformation theory, e.g. of complex structures, in the manner of
Kodaira-Spencer-Grothendieck (cf. e.g. \cite{K2, S} and references
therein), is commonly couched in terms of a differential graded Lie
algebra or dgla $\g$. It can be viewed, as in \cite{cid}, as
studying the \emph{deformation ring} $R(\g)$, defined in term of the
\emph{Jacobi cohomology}, i.e. the cohomology of the Jacobi complex
associated to $\g$. This setting is somewhat restrictive, e.g. it is
not broad enough to accomodate such naturally occurring deformation
problems as embedded deformations of a submanifold $X$ in a fixed
ambient space $Y$. In \cite{atom} we introduced the notion of
\emph{Lie atom} (essentially, Lie pair) and an associated
\emph{Jacobi-Bernoulli complex} as an extension of that of dgla and
its Jacobi complex, one that is broad enough to handle embedded
deformations and a number of other problems besides.\par A purpose
of this paper is to establish the familiar notion of (dg)
\emph{semi-simplicial} Lie algebra (\sela) as an appropriately
general and convenient setting for deformation theory. As a first
approximation, one can think of \sela as a structure like that of
the \v{C}ech complex of a sheaf of Lie algebras on a topological
space $X$ with respect to some open covering of $X$. Not only is
\sela a broad generalization of Lie atom, it is broad enough, as we
show, to encompass deformations of (arbitrarily singular) algebraic
schemes (over $\C$).\par To express the deformation theory of a
\sela $\gb$ we introduce a complex that we call the Jacobi-Bernoulli
complex of $\gb$, though a more proper attribution would be to
Jacobi-Bernoulli-Baker-Campbell-Hausdorff. In a nutshell, the point
of this complex is that it transforms a gluing condition from
nonabelian coycle condition to ordinary (additive) cocycle condition
via the multilinearity of the groups making up the complex. A
typical gluing condition looks like
\begin{equation}\label{nonabelian}
\Psi_{\alpha\beta}\Psi_{\beta\gamma}\Psi_{\gamma\alpha}=1
\end{equation} with $\Psi_{\alpha\beta}\in\exp(\g_{\alpha\beta})$,
where $\g_\ab$ may be thought of as the component of our \sela
$\gb$ having to do with gluing over $U_\alpha\cap U_\beta$. This
condition can be transformed as follows. Write
$$\Psi_{\alpha\beta}=\exp(\psi_{\alpha\beta})$$ etc. Now the BCH
formula gives a formal expression
$$\exp(X)\exp(Y)\exp(Z)=\sum W_{i,j,k}(X,Y,Z)$$ where
$W_{i,j,k}(X,Y,Z)$ is a homogeneous ad-polynomial of tridegree
$i,j,k$ (`BCH polynomial'), which can be viewed as a linear map
$$w_{ijk}:\Sym^i(\g_\ab)\otimes\Sym^j(\g_\bg)\otimes\Sym^k(\g_\ga)\to
g_\abg$$ Then (\ref{nonabelian}) becomes the additive condition that
$$w_{i,j,k}(\psi_\ab^i\otimes\psi_\bg^j\otimes\psi_\ga^k)=0,\ \forall
i,j,k.$$ Now our Jacobi-Bernoulli complex $J(\gb)$ for the \sela
$\gb$ is essentially designed so as to encompass the various BCH
polynomials $w_{ijk}$. It is a comultiplicative complex whose groups
essentially constitute the symmetric algebra on $\gb$ and whose maps
are essentially derived from the $w_{ijk}$ by the requirement of
comultiplicativity. The dual of the cohomology of $J(\gb)$ yields
the deformation ring associated to the \sela $\gb$.\par As mentioned
above, our other main result here is that the deformation theory of
an algebraic scheme over $\C$ can be expressed in terms of a \sela.
Unsurprisingly, this is done via an affine covering. Thus the first
step is to associate a dgla to a closed embedding $$X\to P$$ where
$P$ is an affine (or for that matter, projective) space. We call
this the \emph{tangent} dgla to $X$ and denote it $\t_X(P)$. In a
nutshell, $\t_X(P)$ is defined as the mapping cone of a map that we
construct
$$T_P\otimes\O_X\to N_{X/P}$$ where $N_{X/P}$ is the normal atom to
$X$ in $P$ as in \cite{atom}. That is, $\t_X(P)$ is represented by
the mapping cone of a map of free modules representing
$T_P\otimes\O_X$ and $N_{X/P}$. We will show $\t_X(P)$ admits a dgla
structure, a dgla action on $\O_X$, as well as $\O_X$-module
structure. Up to a certain type of 'weak equivalence', the dgla
$\t_X(P)$ depends only on the isomorphism class of $X$ and not on
the embedding in $P$.
\par The partial independence on the embedding is good enough to
enable us to associate a global \sela $\t_{X\bullet}$ for an
arbitrary algebraic scheme $X$ defined in terms of, but up to weak
equivalence independent of,  an affine covering $X_\alpha$ and
embeddings of each $X_\alpha$ in an affine space $P_\alpha$: e.g.
$$\t_{X,\alpha}=\t_{X_\alpha}(P_\alpha),$$
$$\t_{\X, \ab}=\t_{\X_\alpha\cap X_{\beta}}(P_\alpha\times
P_\beta)$$ etc. Global deformations of $X$ then amount to a
collections of deformations of each $X_\alpha$, given via
Kodaira-Spencer theory by a suitable element
$\phi_\alpha\in\t_{X,\alpha}^1$, plus a collection of gluing data
$\psi_\ab\in\t_{X,\ab}^0,$ and the necessary compatibilities are
readily expressed as a cocycle condition in the Jacobi-Bernoulli
complex $J(\t_{X\bullet}).$

\section{Semi-Simplicial Lie algebras and Jacobi-Bernoulli complex }
\subsection{\sela} Our
notion of \sela is essentially the dual of the portion of the usual
notion of simplicial Lie algebra involving only the face maps
without degeneracy. Let $A$ be a totally ordered index-set. A
\emph{simplex} in $A$ is a finite nonempty subset $S\subset A$,
while a \emph{biplex} is a pair $(S_1\subset S_2)$ of simplices with
$|S_1|+1=|S_2|$; similarly for triplex $(S_1\subset S_2\subset S_3)$
etc. The \emph{sign} $\epsilon(S_1,S_2)$ of a biplex $(S_1,S_2)$ is
defined by the condition that
$$\epsilon((0,...,\hat{p},...,n), (0,...,n))=(-1)^{n-p}.$$ By a
\emph{simplicial Lie algebra} (SLA) $\gb$ on $A$ we shall mean the
assignment for each simplex $S$ on $A$ of a Lie algebra $\g_S$, and
for each biplex $(S_1, S_2)$ of a map ('coface' or 'restriction')
$$r_{S_1,S_2}:\g_{S_1}\to\g_{S_2}$$ such that
$\epsilon(S_1,S_2)r(S_1,S_2)$ is a Lie homomorphism and such that
for each $S_1\subset S_3$ with $|S_1|+2=|S_3|$, we have
\begin{equation}
\label{complex}\sum_{\stackrel{\rm triplex} {(S_1\subset S_2\subset
S_3)}}r_{(S_2,S_3)} r_{(S_1,S_2)}=0.\end{equation} The identity
(\ref{complex}) implies that we may assemble the $\g_S$ into a
complex $K^.(\gb)$ where
$$K^i(\gb)=\bigoplus\limits_{|S|=i+1}\g_S$$ and differential
constructed from the various $r_{(S_1,S_2)}$.
\begin{example}{\begin{rm} If $\g$ is a sheaf of Lie algebras on a
topological space $X$, and $(U_\alpha)$ is an open covering of $X$,
there is a \v{C}ech \sela $$S\mapsto \g(\bigcap\limits_{\alpha\in
S}U_\alpha).$$ The standard complex $K^.(\gb)$ is in this case the
\v{C}ech complex $\check{C}(\g,(U_\alpha))$. This plays a
fundamental role in the study of $\g$-deformations.
\end{rm}}\end{example}
\mex{If $\g\to\h$ is a Lie pair (more generally, a Lie atom, cf.
\cite{atom}), we get a \sela $\g_.$ on $(01)$ with $\g_0=\g, \g_1=0,
\g_{01}=\h.$
$$\h$$
\begin{picture}(80,20)(20,10)
\put(210,30){$\g_\bullet$\line(1,0){30}$_\bullet
0$}\end{picture}\newline
 The deformation-theoretic
significance of $\g_.$ is like that of the Lie atom $(\g,\h)$, viz.
$\g$-deformations together with an $\h$-trivialization.\par An
obvious generalization would be to take a pair of maps $\g_1\to\h,
\g_2\to\h$ (e.g. twice the same map), which corresponds to pairs
($\g_1$- deformation, $\g_2$-deformation) that become equivalent as
$\h$-deformations.}
\subsection{Bernoulli numbers and Baker-Campbell-Hausdorff}   Let $\g$ be a
nilpotent Lie algebra. For an element $X\in\g$ we consider the
formal exponential $\exp(X)$ as an element of the formal enveloping
algbera $\U(\g).$ Then we can write
\begin{equation}\label{beta}\exp(X)\exp(Y)=\exp(\beta(X,Y))\end{equation} where $\beta$ is a certain
bracket-polynomial in $X,Y$, known as the Baker-Campbell-Hausdorff
or BCH polynomial. We denote by $\beta_{i,j}, \beta_i$ the portion
of $\beta$ in bidegree $i,j$ (resp. total degree $i$). Note that
each $\beta_{i,j}(X,Y)$ will be a linear combination  of
(noncommutative) $\ad$ monomials with a total of $i$ many $X$'s and
$j$ many $Y$'s. We write such a monomial in the form
\beq\label{admon1}\ad_S(X^iY^j)=\ad(T_1)\circ...\circ\ad(T_{i+j-1})(T_{i+j})\eeq
where $S\subset [1,i+j]$ is a subset of cardinality $i$ and $T_k=X$
(resp. $T_k=Y$) iff $k\in S$ (resp. $k\not\in S$). We denote by
\beq\ad_S(X_1,...,X_i, Y_1,...,Y_j)\eeq the analogous
function, 
obtained by replacing the $x$th occurrence of $X$ (resp. $y$th
occurrence of $Y$) by an $X_x$ (resp. $Y_y$) and by
$\ad_S^\Sym(X_1,...,X_i,Y_1,...,Y_j)$ the corresponding symmetrized
version, i.e.
\beq\label{admon3}\ad_S^\Sym(X_1,...,X_i,Y_1,...,Y_j)=\sum\limits_{
\stackrel{\pi\in\frak S_i}{\rho\in\frak S_j}}\frac{1}{i!j!}
\ad_S(X_{\pi(1)},...,X_{\pi(i)},Y_{\rho(1)},...,Y_{\rho(j)}).\eeq\par
We will compute $\beta$, following \cite{vara}, \S2.15 (where
Varadarajan attributes the argument to lectures of Bargmann that
follow original papers by Baker and Hausdorff). Set
$$D(x)=\frac{e^x-1}{x}, C(x)=1/D(x).$$ Thus, $C(x)$ is the
generating function for the Bernoulli numbers $B_n$, i.e.
$$C(x)=1+\sum\limits_{n=1}^\infty \frac{B_n}{n!}x^n
=\sum\limits_{n=0}^\infty C_nx^n.$$ Now the reader can easily check
that for any derivation $\partial$ we have
$$\partial\exp(U)\exp(-U)=D(\ad (U))(\partial U),
\exp(-U)\partial\exp(U)=D(-\ad(U))(\partial U).$$ Now differentiate
(\ref{beta}) with respect to $X$ and multiply both sides by
$\exp(-\beta(X,Y)).$ This yields (where $\partial_X$ is the unique
derivation taking $X$ to $X$ and $Y$ to 0)\par 
\lefteqn{ X=\partial_X(\exp(X))\exp(-X)=}
$$~~~~~~~~~~~~~~~~~~
{=\partial_X(\exp(\beta(X,Y)))\exp(-\beta(X,Y))
=D(\ad(\beta(X,Y)))(\partial_X\beta(X,Y)).}$$
Thus
\begin{equation}\label{betax}
\partial_X\beta(X,Y)=C(\ad(\beta(X,Y)))(X).
\end{equation} Similarly,
\begin{equation}\label{betay}\partial_Y\beta(X,Y)=
C(-\ad(\beta(X,Y)))(Y).
\end{equation} Starting from $\beta_0=0$, the formulas
(\ref{betax}),(\ref{betay}) clearly determine $\beta$. For example,
clearly $\beta_{0,*}(X,Y)=Y$, therefore it follows that
\begin{equation}\beta_{1,*}(X,Y)=
C(\ad(Y))(X)=X+\frac{1}{2}[X,Y]+\frac{1}{12}\ad(Y)^2(X)+...
\end{equation} \par We
shall require the obvious extension of this set-up to the trivariate
case. Thus define a function $\beta(X,Y,Z)$ (NB same letter as for
the bivariate version) by
\beq\exp(X)\exp(Y)\exp(Z)=\exp(\beta(X,Y,Z))\eeq and let
$\beta_{i,j,k}$ denote its portion in tridegree $(i,j,k)$. Note that
$$\beta(X,Y,Z)=\beta(\beta(X,Y),Z).$$
\subsection{Jacobi-Bernoulli complex}
Let $\gb$ be a \sela. For simplicity, we shall assume $\gb$ is
2-dimensional, in the sense that $\g_S=0$ for any simplex $S$ of
dimension $>2$; for our applications to deformation theory, this is
not a significant restriction. We will also assume that $\gb$ is
\emph{strongly nilpotent} in the sense that it is an algebra over a
commutative ring $R$ such that $\g_S^{\otimes N}=0$ for all
simplices $S$ and some integer $N$ independent of $S$, with all
tensor products over $R$. This condition obviously depends only on
the $S$-module structure of $\gb$ and not on its Lie bracket. We are
going to define a filtered complex $J=\jsh_m(\gb).$ The groups $J^j$
can be defined succinctly as
$$J^j=(\Sym^.(K^.(\gb)[1]))^j$$ where $\Sym^.$ is understood in the
signed or graded sense, alternating on odd terms, and $K^.(\gb)[1]$
is the standard complex on $\gb$ shifted left once (which is a
complex in degrees $-1,0,1$) . The increasing filtration $F.$ is by
`number of multiplicands'.  More concretely,\beq
J^{j,k}=\bigoplus\limits_
{\stackrel{-\sum\limits_i\l_i+\sum\limits_in_i=j}
{\sum\limits_i\l_i+\sum\limits_im_i+\sum\limits_in_i=
k}}\bigotimes\limits^i\bigwedge\limits^{\l_i}\g_{\alpha_i} \otimes
\bigotimes\limits^i\Sym^{m_i}\g_{\alpha_i\beta_i}\otimes
\bigotimes\limits^i\bigwedge\limits^{n_i}\g_{\alpha_i\beta_i\gamma_i}
\eeq \beq F_mJ^j=\bigoplus\limits_{k\leq m}J^{j,k},\eeq \beq
J^j=F_\infty J^j=F_NJ^j.\eeq To define the differential $d$ on
$J^.$, we proceed in steps. Let $\alpha<\beta<\gamma$ be indices and
recall that we are identifying $\g_{\gamma\alpha}$ with
$\g_{\alpha\gamma}$.\begin{itemize}\item The differential is defined
so that the obvious inclusion \beq K^.(\gb)[1])=F_1J^.\to J^.\eeq is
a map of complexes.
\item The
component\begin{equation*} \Sym^i\g_{\gamma\alpha}\otimes
\Sym^j\g_{\ab}\otimes\Sym^\g_{\beta\gamma}\to
\g_{\alpha\beta\gamma}\end{equation*} is given by\begin{equation}
 {X^iY^jZ^k\mapsto\beta_{i,j,k}(X,Y,Z)}.\end{equation}
\item The component
\begin{equation*}\nonumber\g_\alpha\otimes\Sym^i\g_{\alpha\beta}\otimes
\Sym^n\g_{\beta\gamma}
\to\Sym^{i-t+1}\g_{\alpha\beta}\otimes\Sym^n\g_{\beta\gamma}, 0\leq
t\leq i\end{equation*} is given by\begin{equation} X\otimes
Y^i\otimes Z^n\mapsto C_tY^{i-t}\ad(Y)^t(X)\otimes
Z^n\end{equation}(where $C_t$ is the normalized Bernoulli
coefficient).\item Other componets are defined subject to the
'derivation rule', e.g. the component
$$\g_\alpha\otimes\Sym^i\g_{\beta\gamma}\otimes\Sym^k\g_{\beta\gamma}
\to\g_{\gamma\alpha}\otimes
\Sym^i\g_{\ab}\otimes\Sym^k\g_{\beta\gamma}$$ is extended in the
obvious way from the given differential
$\g_\alpha\to\g_{\gamma\alpha}$.\item Components not defined via the
above rules are set equal to 0. In particular, the component
$$\g_\alpha\otimes\g_{\alpha\beta\gamma}\to \g_{\alpha\beta\gamma}$$
is zero.
\end{itemize} The following result summarizes the main properties of
the Jacobi-Bernoulli complex $J$ associated to a \sela (not least,
that it is a complex!). It is in part, but not entirely, a direct
extension of the analogous result for Lie atoms given in
\cite{atom}.
\begin{thm}\label{tang sela}
\begin{enumerate}\item $(J^., F.)$ is a functor from the
category of SELAs over $S$ to that of comultiplicative,
cocummutative and coassociative filtered complexes over $S$.\item
The filtration $F_.$ is compatible with the comultiplication and has
associated graded
$$F_i/F_{i-1}=\bigwedge\limits^i(\gb).$$\item
 The quasi-isomorphism class of $J(\gb)$ depends only on the
quasi-isomorphism class of $\gb$ as \sela.\end{enumerate}
\end{thm}\begin{proof} As in the proof of \cite{atom}, Thm 1.2.1,
the main issue is to prove $J$ is a complex, i.e. $d^2=0$. And again
as in \cite{atom}, it suffices, in light of the derivation rule,  to
prove the vanishing of the components of $d^2$ that land in $F_1$,
i.e. that have just one multiplicative factor. Among those, the
proof that these components of $d^2$ vanish on terms of degree $\leq
-2$, i.e. involving $\bigwedge\limits^i\g_\alpha, i\geq 2$, is
similar to the case of the JB complex considered in \cite{atom}. The
essential new case, not considered in \cite{atom}, is the vanishing
of the $F_1$- components of $d^2$ on terms of degree $-1$, i.e.
terms of the from
$$X\otimes Y^i\otimes Z^n\in
\g_\alpha\otimes\Sym^i\g_{\alpha\beta}\otimes
\Sym^n\g_{\beta\gamma}.$$ For such a term, what needs to be shown is
the vanishing of the component of $d^2$ of it in
$\g_{\alpha\beta\gamma}$. Thus, we need to prove that\beq d^2(
X\otimes Y^i\otimes Z^n)_{\g_{\alpha\beta\gamma}}=0.\eeq Now this
component gets contributions via the various components of $d(
X\otimes Y^i\otimes Z^n)$ and those contributions come in  two
kinds:\begin{itemize}\item Via
$g_{\gamma\alpha}\otimes\Sym^i\g_{\alpha\beta}\otimes
\Sym^n\g_{\beta\gamma}$, we get $-\beta_{1,i,n}(X,Y,Z)$. This comes
from $\beta(\beta(X,Y),Z)$, but is only affected by the terms in
$\beta(X,Y)$ of degree $\leq 1$ in $X$, i.e. by
$$U=Y+C(\ad(Y))(X)=Y+\sum\limits_{t=0}^\infty C_t\ad(Y)^t(X).$$ This
contribution is obtained by taking $-\beta_{i+1-t,n}(U,Z)$ and
replacing each monomial $$\ad_S(U^{i+1-t}Z^n)$$ (cf.
(\ref{admon1}-\ref{admon3}) ) by
$$(i+1-t)\ad_S^\Sym(C_t(\ad(Y)^t(X))Y^{i-t}Z^n)$$
and finally summing as $t$ ranges from 0 to $i$.
\item Via
$\Sym^{i+1-t}\g_{\alpha\beta}\otimes\Sym^n\g_{\beta\gamma}$, for
each $0\leq t\leq i$, we get a contribution equal to the expression
obtained by taking $\beta_{i+1-t,n}(W,Z)$, replacing each monomial
$\ad_S(W^{i+1-t}Z^n)$ by
$(i+1-t)\ad_S^\Sym(C_t\ad(Y)^t(X)Y^{i-t}Z^n)$.
\end{itemize} Thus, the sum total of all contributions to
$d^2(X\otimes Y^i\otimes Z^n)_{\g_{\alpha\beta\gamma}}$ is zero.
\end{proof} The ring
$$R(\gb)=\C\oplus\HH^0(J(\gb))^*$$ is called the \emph{deformation
ring} of $\gb$.
\subsection{Special multiplicative cocycles} Let $(S,\m_S)$ be a
local artin $\C$-algebra and let $\gb^.=\gb$ be a dg-\sela (i.e.
each $\g_T$ is a dgla and the coface maps are dg homomorphisms). A
special class of (multiplicative) cocycles for the Jacobi-Bernoulli
complex $$J^.(\gb\otimes\m_S)\subset J^.(\gb)\otimes\m_S$$ can be
constructed as follows. Suppose $$\phi_\bullet\in
K^0(\gb)^1\otimes\m_S=\bigoplus\limits_\rho\g_\rho^1\otimes\m_S,$$
$$\psi_\bullet\in K^1(\gb)^0\otimes\m_S=\bigoplus_\rs\g_\rs^0\otimes\m_S$$
are such that, $\forall \rho,\sigma,\tau$,
$$\partial\phi_\rho=-\half [\phi_\rho,\phi_\rho],$$
$$\partial\psi_\rs=C(\ad(\psi_\rs))(\phi_\rho-\phi_\sigma),$$
$$\beta(\psi_\rs, \psi_{\sigma\tau}, \psi_{\tau\rho})=0,\
 .$$
Then let $$\epsilon(\phi_\bullet, \psi_\bullet)\in
J^0(\gb)\otimes\m_S$$ be the element with components $$\psi_\rs\in
\g_\rs\otimes\m_S,$$$$\phi_\rho\in\g_\rho\otimes\m_S, $$ and
generally $$\bigwedge\limits^r\phi_\rho\otimes(\psi_\rs)^n\in
\bigwedge\limits^r(\g^1_\rho\otimes\m_S)
\otimes\Sym^n(\g_\rs\otimes\m_S), r, n \geq 0.$$ We call
$\epsilon(\phi_\bullet, \phi_\bullet)$ a \emph{special
multiplicative cocycle with coefficients in $S$}.
\begin{lem}\begin{enumerate}\item
 The cochain $\epsilon(\phi_\bullet, \psi_\bullet)$
defined above is a 0-cocycle for $J(\gb)$ and the associated map
$$^t\epsilon(\phi_\bullet, \psi_\bullet)\in\mathrm{Hom}( R(\gb),S)$$
is a local  ring homomorphism.\item Given $S$, there is a bijection
between cohomology classes of special multiplicative cocycles with
coefficients in $S$ and local ring homorphisms $R(\gb)\to
S.$\end{enumerate}
\end{lem}\qed

\section{Tangent \sela}\subsection{Affine schemes: tangent
dgla }We work in the 'global affine' setting, though the
construction can evidently be coherently sheafified. Let $X$ be a
closed subscheme of an affine space $P$ and let $I=I_{X/P}$ denote
the ideal of $X$ in the coordinate ring $A_P$ (a similar
construction can be done for $P$ an arbitrary open subscheme of a
projective space). Let\begin{equation}...\to F^{-1}\to F^0\to
I\end{equation} be a free resolution of $I$. Thus each $F^i$ is a
free module on generators $e_\alpha^i$ which correspond for $i=0$ to
generators $f_\alpha$ of $I$ and for $i<0$ to syzygies $\sum
f_{\alpha\beta}e^{i+1}_\beta$.
\par
 Next, set $F^1= A_P$ (the coordinate ring) and let $F^._+$ be the complex in
degrees $\leq 1$ given by
$$...\to F^0\stackrel{\epsilon}{\to} F^1$$ which is a free
resolution of $A_X$. Then set \begin{equation}
\n=\n_{X/P}=\hom^.(F^., F^._+)\end{equation} which we view as
sub-dgla of $\endo^.(F^._+)$ consisting of maps vanishing on
$F^1_+$. This may be called the \emph{normal dgla} of $X$ in $P$.
Next, we define a map $$\kappa: T_P\to \n^1=\bigoplus\limits_{i\leq
0} \hom(F^i_+\to F^{i+1}_+)$$ as follows:\begin{itemize}\item For
$v\in T_P$, $\kappa(v)^0$ is the map taking a distinguished
generator $e_\alpha^0$ of $F^0$, which corresponds to a generator
$f_\alpha$ of $I$, to $v(f_\alpha)$.\item For $i<0$, $\kappa(v)^i$
takes a distinguished generator $e_\alpha^i$ of $F^i$, which
corresponds to a syzygy $\sum f_{\alpha\beta}e^{i+1}_\beta$, to
$\sum v(f_{\alpha\beta})e^{i+1}_\beta$
\end{itemize}
It is not hard to see that this yields a map of complexes
$$T_P\to\n[1],$$ i.e. that $d_\n\circ\kappa=0$.
 Then proceeding similarly, one can lift $\kappa$
to a map of complexes \begin{equation}\kappa: T_P\otimes
F^._+\to\n[1]\end{equation} The mapping cone carries a natural
structure of dgla which we call the \emph{tangent dgla} of $X$ with
reference to $P$ and denote by $\t_X(P)$. It is not hard to see that
$\t_X(P)$ is acyclic in negative degrees and bounded. It admits an
'$\O_X$-module structure' in the form of a pairing $$F^._+\times
\t_X(P)\to\t_X(P),$$ as well as a 'derivation action' on $\O_X$
(i.e. on $F^._+$).\par As for the dependence on the embedding, let
$X\to Q$ be another affine embedding. Then via the diagonal we get a
third one $X\to P\times Q$, together with dgla maps
$$\t_X(P)\to\t_X(P\times Q), \t_X(Q)\to \t_X(P\times Q)$$ which, it
is easy to check, induce isomorphisms on cohomology in degree $\leq
1$ and an injection on $H^2$. We call a dgla morphism with these
properties a \emph{direct weak equivalence} and define a general
\emph{weak equivalence} of dgla's  to be a composition of weak
equivalences and their inverses.\begin{example}{\rm If $X$ is a
hypersurface with equation $f$ in $P=\A^n$, its tangent dgla may be
identified with the complex in degrees $-1,0,1$
$$nA_P\stackrel{(f,0)}{\to}nA_P\oplus A_P\stackrel{(\partial f/\partial
x_1,...,\partial f/\partial x_n, f)}{\longrightarrow}A_P$$ Its $H^1$
is the so-called Milnor algebra of $f$ (finite-dimensional if $f$
has isolated singularities). }\end{example}\subsection{Maps of
affine schemes: tangent dgla} The notion of tangent dgla of an
affine scheme can be extended to the case of a \emph{mapping} of
affine schemes, as follows. Let
$$f:X\to Y$$ be a mapping of affine schemes. Given affine embeddings
$X\to P, Y\to Q$, $f$ can be extended to a map $P\to Q$. Replacing
$X\to P$ by the graph embedding $X\to P\times Q$, we may assume
$P\to Q$ is a product projection. Then we have an injection
$I_{Y,Q}\to I_{X,P}$ which extends to the free resolutions $F^._Y\to
F^._X$, and we may moreover assume that each $F^i_Y\otimes A_P\to
F^i_X$ is a direct summand inclusion. We can identify the functor
$f^!$ on complexes with $f^!\cdot =\cdot\otimes_{A_Y} F^._{X+}.$
Then the complex $f^!(\n_{Y/Q})$ can be represented by
$$\hom^._{A_Q}(F^._Y, F^._{+X})=\hom_{A_P}^.(F^._Y\otimes A_P,
F^._{+X})$$ and there are maps \begin{equation}\label{nf}\n_{X/P}\to
f^!(\n_{Y/Q})\leftarrow \n_{Y/Q}\end{equation} The mapping cone of
\ref{nf} can be represented by the sub-dgla of
$\n_{X/P}\oplus\n_{Y/Q}$ consisting of pairs $(a^.,b^.)$ such that
$a^.$ vanishes on the subcomplex $F^._Y\otimes A_P\subset F^._X$. We
denote this mapping cone by $\n_f$ or more properly $\n_{f,P,Q}$ and
refer to it as the \emph{normal dgla} of $f$.\par Next, proceeding
as in the case of schemes, we can construct a suitable
representative of the mapping cone $K$ of
$$T_P\otimes F^._{+X}\to T_Q\otimes F^._{+X}\leftarrow T_Q\otimes
F^._{+Y},$$ together with a map of $K$ to $\n_f$, so that the
mapping cone of $K\to\n_f$ is a dgla, called the \emph{tangent dgla}
to $f$ and denoted $\t_f$ or more properly, $\t_f(P,Q).$ By
construction, $\t_f(P,Q)$ is the mapping cone of
\begin{equation}\label{tfseq}\t_X(P)\oplus\t_Y(Q)\to
f^!\t_Y(Q).\end{equation}
\par
\subsection{Schemes and maps: tangent \sela}
 Here we construct the tangent \sela of a
separated algebraic scheme $X$ over $\C$ (the separatedness does not
seem to be essential). Let $(X_\rho)$ be an affine open covering of
$X$ indexed by a well-ordered set, and for each $\rho$ let $P_\rho$
be an affine space with a closed embedding
\begin{equation}\iota_\rho:X_\rho\subset P_\rho.\end{equation}
 Set $B_\rho=A_{P_\rho}$. We call the
system $(X_\rho\subset P_\rho)$ an \emph{affine embedding system}
for $X$. Via the diagonal, $X_\rho\cap X_\sigma$ is a closed
subscheme of $X_\rho\times X_\sigma$, hence of $ P_\rho\times
P_\sigma.$ Similarly, for any multi-index $\rho_0<...<\rho_k$, we
define
\begin{equation} X_{(\rho_0,...,\rho_k)}=
X_{\rho_0}\cap...\cap X_{\rho_k},\ \
P_{(\rho_0,...,\rho_k)}=P_{\rho_0}\times...\times P_{\rho_k}
\end{equation}
and the natural closed embedding
\begin{equation}\iota_{(\rho_0,...,\rho_k)}:\label{k-fold inters}
X_{(\rho_0,...,\rho_k)}\subset
P_{(\rho_0,...,\rho_k)},\end{equation} and we denote the ideal of
the latter by $I_{(\rho_0,...,\rho_k)}$. We call the system \beq
\large{(} X_{(\rho_0,...,\rho_k)}\subset P_{(\rho_0,...,\rho_k)},
(\rho_0<...<\rho_k), k\geq 0\large{)}\eeq the \emph{simplicial
extension} of the affine embedding system $(X_\rho\subset P_\rho)$.
Note that the defining equations for the image of
$\iota_{(\rho_0,...,\rho_k)}$ consist of defining equations for the
images of individual embeddings $\iota_{\rho_i}$, together with
equations for the small diagonal on $X^{k+1}.$ The latter are of
course generated by the pullbacks of the equations of the small
diagonal in $X^k$ via the various coordinate projections $X^{k+1}\to
X^k$. Therefore, it is possible to choose mutually compatible free
resolutions for all the $I_{(\rho_0,...,\rho_k)}$, and we denote
these by $F^._{(\rho_0,...,\rho_k)}$. In fact, we may assume that
\beq
F^1_{(\rho_0,...,\rho_k)}=B_{(\rho_0,...,\rho_k)}:=\bigotimes\limits_0^k
B_{\rho_i},\eeq \beq
F^i_{(\rho_0,...,\rho_k)}=\bigoplus\limits_{j=0}^k
(F^i_{\rho_j}\otimes
B_{(\rho_0,...,\rho_k)})\oplus\Delta^i_{(\rho_0,...,\rho_k)}, i\leq
0,\eeq where $\Delta^._{(\rho_0,...,\rho_k)}$ is a lifting to
$B_{(\rho_0,...,\rho_k)}$ of a free resolution of the small
diagonal\nl
$X_{(\rho_0,...,\rho_k)}\subset\prod\limits_{j=0}^kX_{\rho_j}$ and
moreover for any biplex
$$\rho^k=(\rho_0,..,\rho_k)\subset \rho^{k+1}=
(\rho_0,...,\rho_{k+1}),$$ if we let
$$\pi_{\rho^{k+1},\rho^{k}}:P_{\rho^{k+1}}\to P_{\rho^{k}}$$ denote
the natural projection, then we have a direct summand inclusion
\begin{equation}\label{Odiff}\pi_{\rho^{k+1},\rho^{k}}^*F^._{\rho^{k}}:=
\pi_{\rho^{k+1},\rho^{k}}\inv F^._{\rho^{k}}\otimes
B_{{\rho^{k+1}}}\to F^._{\rho^{k+1}}\end{equation} Putting together
these groups and maps, and twisting by the appropriate sign, i.e.
$\epsilon(\rho^k,\rho^{k+1})$, we get an 'extrinsic \v{C}ech
(double) complex'
\begin{equation} \xc(\O_X):
\bigoplus\limits_{\rho^0}F^._{+\rho^0}\to...\to
\bigoplus\limits_{\rho^k}F^._{+\rho^k}\to...\end{equation} Actually
this is quasi-isomorphic to the usual \v{C}ech complex of $\O_X$,
but we can do more with it. Note that the map \ref{Odiff} gives rise
to a direct summand inclusion
$$\pi_{\rho^{k+1},\rho^{k}}^*(\gl^.(I_{\rho^{k}}))\to
\gl^.(I_{\rho^{k+1}}),$$ whence a dgla map
$$\delta_{\rho^{k},\rho^{k+1}}:
\n_{X_{\rho^{k}}/ P_{\rho^{k}}}\to
\n_{X_{\rho^{k+1}}/P_{\rho^{k+1}}}$$ Then we can similarly construct
a 'normal \sela ' $$\n_{X_\bullet/P_\bullet}:...\to
\bigoplus\limits_{\rho^k}\n_{X_{\rho^{k}}/ P_{\rho^{k}}}\to...$$
Likewise, we have an 'ambient tangent complex' $T_{P_\bullet}$ and
$T_{P_\bullet}\otimes\xc(\O_X)$ and a map
\begin{equation}\label{tangsela}
T_{P_\bullet}\otimes\xc(\O_X)\to
\n_{X_\bullet/P_\bullet}\end{equation} Finally, we define the
\emph{tangent \sela} of $X$ (with reference to the simplicial system
$(X_\bullet, P_\bullet)$) to be the mapping cone of this, and denote
it by $\t_{X\bullet}(X_\bullet, P_\bullet)$ or simply
$\t_{X\bullet}$. This is the \sela whose value on the simplex
$\rho^k$ is the dgla $\t_{X_{\rho^k}}(P_{\rho^k})$. By Theorem
\ref{tang sela}, there is an associated Jacobi-Bernoulli complex
$J^.(\t_{X_\bullet})$, which we denote by $J^._X$ and refer to as
the \emph{Jacobi-Bernoulli complex of} $X$. Up to filtered,
comultiplicative quasi-isomorphism, it depends only on the
isomorphism class of $X$ as scheme over $\C$. Therefore the
\emph{deformation ring} of $X$ $$R_X=\C\oplus\HH^0(J^._X)^*$$ is
canonically defined. In the next section we relate $R_X$ to flat
deformations of $X$ over artin rings. For any artin local
$\C$-algebra $S$, we set
$$J^._{X,S}=J^.(\t_{X_\bullet}\otimes\m_S)$$ and note that via the
natural map $J^._{X,S}\to J^._X\otimes\m_S$, any class
$\epsilon\in\HH^0(J^._{X,S})$ yields a local homomorphism
('classifying map')
$$^t\epsilon: R_X\to S.$$
\par As in the affine case, this construction may be extended
to the case of maps. Thus let
$$f:X\to Y$$ be a morphism of schemes. Then we can choose respective
affine coverings
$$X_\alpha\to P_\alpha, Y_\alpha\to Q_\alpha \textrm{\ \ such that\
\ } f(X_\alpha)\subset Y_\alpha.$$ Then for each simplex $\rho^k$,
the restriction of $f$ yields a morphism
$$f_{\rho^k}:X_{\rho^k}\to Y_{\rho^k},$$
and for this we have an associated tangent dgla $\t_{f_{\rho^k}}$.
Putting these together, we get a tangent \sela (with respect to the
given affine coverings)
$$\t_f:...\to \t_{f_{\rho^k}}\to...$$ As before, $\t_f$ is the
mapping cone of $$\t_{X\bullet}\oplus\t_{Y\bullet}\to
f^!\t_{Y\bullet}.$$ Thus we have \sela morphisms $$\t_f\to \t_Y,
\t_f\to \t_X, f^!\t_Y[1]\to\t_f.$$ Correspondingly, we have a
Jacobi-Bernoulli complex $J^._f$, a deformation ring $R_f$ together
woth maps $R_X\to R_f, R_Y\to R_f.$\begin{rem}{\rm When $X$ is
smooth, its tangent \sela is equivalent to a dgla, e.g. the
Kodaira-Spencer algebra, a soft dgla resolution of the tangent
sheaf.}\end{rem}
\par\section{Deformations of schemes} \subsection{Classification}
Deformations of an algebraic
scheme $X/\C$ can be classified in terms of the associated tangent
\sela $\t_X$ and its Jacobi-Bernoulli cohomology. Consider first the
case of an affine scheme $X\subset P$ (notations as in \S2.1). Let
$S$ be a local artinian $\C$-algebra. Then a flat deformation of $X$
over $S$ is determined by, and determines, up to certain choices, an
element
\begin{equation}\label{aff_ks_elt}\phi\in\t^1_X(P)\otimes\m_S=\hom^1(F^._X,
F^._{+X})\otimes\m_S\end{equation} known as a \emph{Kodaira-Spencer
cochain}, which satisfies the integrability condition
\begin{equation}\label{aff_integ}\partial\phi=-\half[\phi,\phi].\end{equation}
The deformation corresponding to $\phi$ can be determined e.g. as
the subscheme of $P\times\Spec(S)$ having $(F^._X\otimes S,
\partial +\phi)$ as resolution; we may denote this by $X^\phi$.\par
Now globally, let $X$ be an algebraic scheme over $\C$ and as in
\S2.3 choose an affine embedding system $$\iota_\rho:X_\rho\to
P_\rho.$$ This gives rise as in \S2.3 to a representative for the
tangent \sela $\t_X$. Now suppose given a deformation of $X$ over
$S$ as above. This restricts for each $\rho$ to a deformation of
$X_\rho$, whence a Kodaira-Spencer cochain
$\phi_\rho\in\t^1_{X_\rho}(P_\rho)\subset \t^1_X(P_\bullet)$,
satisfying an integrability condition as in (\ref{aff_integ}), so
that the restricted deformation of $X_\rho$ is $X_\rho^{\phi_\rho}$.
Moreover, the fact that $\phi_\rho$ and $\phi_\sigma$ restrict to
equivalent deformations of $X_{\rho\sigma}\subset P_{\rho\sigma}$
yields an isomorphism
\begin{equation}\label{sigma_iso_rho}X_\rho^{\phi_\rho}\cap
X_\sigma\simeq X_\sigma^{\phi_\sigma}\cap X_\rho;\end{equation} both
of these are closed subschemes of $P_{\rho\sigma}\times\Spec(S)$ and
the isomorphism (\ref{sigma_iso_rho}) extends to an automorphism of
$P_{\rho\sigma}\times\Spec(S)$, necessarily of the form
$$\exp(t_{\rho\sigma}), t_{\rho\sigma}\in
T_{P_{\rho\sigma}}\otimes\m_S.$$ Then we get two resolutions of
$X_\rho^{\phi_\rho}\cap X_\sigma$, the 'original' one with
differential $\partial +\phi_\rho$, and the one pulled back from
$X_\sigma^{\phi_\sigma}\cap X_\rho$, whose differential is $\partial
+\phi_\sigma.$ It is easy to see and well known that the two
resolutions differ by an isomorphism of the form $\exp(u_\rs)$ where
$u_\rs\in\gl^0(F^._{X_\rho})\otimes\m_S.$ Thus all in all there is a
(uniquely determined) element
$$\psi_{\rho\sigma}\in\t^0_{X_{\rho\sigma}}\otimes\m_S=
(T_{P_\rs}\oplus\gl^0(F^._{X_\rs}))\otimes\m_S$$ such that
\begin{equation}\label{phipsi_cocycle}
\exp(\psi_{\rho\sigma})(\partial+\phi_\sigma)
\exp(-\psi_{\rho\sigma})=\partial+\phi_\rho.\end{equation} By
construction, we clearly have
\begin{equation}\label{psi_cocycle}
\exp(\psi_{\rho\sigma})\exp(\psi_{\sigma\tau})=
\exp(\psi_{\rho\tau})\end{equation}
Thus, using (\ref{aff_integ}, \ref{phipsi_cocycle},
\ref{psi_cocycle}),
 $\epsilon(\phi_\bullet, \psi_\bullet)$ is a special multiplicative
cocycle in the Jacobi-Bernoulli complex $J(\t_X)\otimes \m_S.$
Conversely, given a special multiplicative cocycle
$\epsilon(\phi_\bullet, \psi_\bullet)$ with values in $S$, the
$\phi_\bullet$ data yields a collection of deformations of the
affine pieces of $X$, while the $\psi_\bullet$ glues these
deformations together. These processes are inverse to each other up
to an automorphism, and are precise mutual inverses when there are
no automorphisms. Hence \begin{thm} Let $X$ be an algebraic scheme
over $\C$ such that $H^0(\t_X)=0.$ Then for any local artin
$\C$-algebra $S$, there is a bijection between the set of
equivalence classes of flat deformations of $X$ over $S$ and the set
of local homomorphisms from $R_X$ to $S$.\end{thm}
\subsection{Obstructions} Let $S$ be a local artin algebra and $I<S$
an ideal contained in the socle $\mathrm{ann}_S(\m_S)$ and
$\bar{S}=S/I.$ Let $\bar{\epsilon}=\epsilon(\bar{\phi}_\bullet,
\bar{\psi}_\bullet)$ be a special multiplicative cocycle with
coefficients in $\bar{S}$. Let $\phi_\bullet, \psi_\bullet$ be
arbitrary liftings of $\phi_\bullet, \psi_\bullet$ with coefficients
in $S$. Thus, $\epsilon=\epsilon(\phi_\bullet, \psi_\bullet)$ is not
necessarily a cocycle. However, it is easy to check that the
coboundary $\partial \epsilon$ lies in  $$(F_1J_X)^1\otimes
I=(K^0(\t_X)^2\oplus K^1(\t_X)^1\oplus K^2(\t_X)^0)\otimes I$$
(because it dies $\mod I$) and moreover, that $\partial\epsilon$ is
a cocycle for $\mathrm{tot}(K^.(\t_X)^.)\otimes I$ (because it is a
cocycle for $J^._{X,S}$). Thus, we obtain a cohomology class
\begin{equation}\label{obstruction} \mathrm{ob}(\bar{\phi}_\bullet,
\bar{\psi}_\bullet)\in H^2(\t_X)\otimes
I=H^2(\mathrm{tot}(K^.(\t_{X})^.)\otimes I.\end{equation} This class
is independent of choices and represents the obstruction to lifting
$\bar{\epsilon}$ to a special multiplicative cocycle with
coefficients in $S$.
\subsection{Applications: relative obstructions, stable
subschemes and surjections} Let $f:X\to Y$ be an embedding of a
closed subscheme. We consider the question of 'relative
obstructions' i.e. obstructions to lifting a given deformation of
$Y$ to a deformation of $f$. These are obstructions for the mapping
cone of $\t_f\to \t_Y$, which is the same as that of $\t_X\to
f^!\t_Y.$ Locally, if
 $Y\to P$ is an affine embedding, such obstructions have values in
$\mathrm{Ext}^1(I_{X/Y}, \O_X)$, where $I_{X/Y}=I_{X/P}/I_{Y/P}$ is
\emph{independent of $P$} even as complex up to quasi-isomorphism
(not just weak equivalence). Hence, global obstructions also have
values in $\mathrm{Ext}^1(I_{X/Y}, \O_X)$. In reasonably good cases
though, the obstruction group can be narrowed considerably. The
following result sharpens Thm. 1.1 of \cite{Rhol}\begin{thm} Let
$X\subset Y$ be a closed subscheme having no component contained in
the singular locus of $Y$. Then obstructions to deforming $X\to Y$
relative to deforming $Y$ are in
\begin{equation}\label{image}\mathrm{Im(Ext}^1_Y(I_{X/Y}/I_{X/Y}^2,
\O_X)\to\mathrm{Ext}^1_Y(I_{X/Y}, \O_X)).\end{equation}\end{thm} It
follows easily, in particular, that if $X\to Y$ is moreover a
regular embedding with normal bundle $N$ and $H^1(X,N)=0$, then
$X\to Y$ is 'relatively unobstructed' or 'stable' relative to $Y$,
i.e. deforms with every deformation of $Y$, and furthermore the
Hilbert scheme of $Y$ is smooth at the point corresponding to $X$.
This generalizes a result of Kodaira \cite{K1} in the smooth
case.\par To sketch the proof, working locally, let $F^._Y$ be a
free resolution of $I_{Y/P}$ and extend it to a free resolution
$F^._X$ of $I_{X/P}$, such that, termwise,
$$F^i_X=F^i_Y\oplus F^i_{X/Y}$$ where $F^i_{X/Y}$ is a suitable free
complement, $F^._Y\to F^._X$ is a map of complexes (though not a
direct summand inclusion), and $F^._{X/Y}\otimes\O_Y$, which is a
quotient complex of $F^._X$, is a free resolution of $I_{X/Y}$. We
may also assume that $F^._{X/Y}$ contains a subcomplex $F^._{X^2}$
(with termwise direct summands) resolving $I_X^2$. A deformation of
$Y$ yields a linear map
$$v:F^0_Y\to \O_X.$$ The obstruction to lifting this to a
deformation of $X$ is given by $$v\circ\delta:F^{-1}_{X/Y}\to \O_X$$
where $\delta:F^{-1}_{X/Y}\to F^0_Y$ is the 'connecting map' from
the resolution. $\delta$ takes a relation among generators of
$I_X\mod I_Y$ to the appropriate linear combination of generators of
$I_Y$. Now, and this is the point, our assumption about
singularities means that no generator of $I_Y$ can be in $I_X^2$, or
more precisely, that
$$I_{X}^2\cap I_Y\subset I_XI_Y.$$ This implies that
$\delta(F^{-1}_{X^2})\subset I_XF^0_Y$. Since $v$ is $\O_P$-linear,
it follows that the composition $v\circ\delta$ is zero on
$F^{-1}_{X^2}$ and lives in the image as in (\ref{image}). This
proves the result in the affine case, and the extension to the
global case is straightforward.\qed.\par Next we consider an
application to surjections is (compare \cite{Rhol}, Thm
2.1):\begin{thm} Let $f:X\to Y$ be a projective morphism with
$$f_*(\O_X)=\O_Y, R^1f_*(\O_X)=0.$$ Then $f$ deforms with every
deformation of $X$.\end{thm}\begin{proof} Here we use the fact that
the mapping cone of $\t_f\to\t_X$ is equivalent to that of $\t_Y\to
f^!\t_Y.$ To prove the result it suffices to show the
following\newline (*)\emph{ the natural map $\t_Y\to f^!\t_Y$
induces a surjection on $H^1$ and an injection on $H^2.$}\par Indeed
the $H^1$-surjectivity property implies that any first-order
deformation of $X$ lifts to a deformation of $f$; then the
$H^2$-injectivity property allows us to extend this inductively, via
obstruction theory (\S3.2), to $n$th order deformations. To prove
(*), note that $f^!\t_Y$ can be represented by the tensor product
$\t_Y\otimes\xc(\O_X)$, which is a double complex with terms
$\t^i_Y\otimes \xc^j(\O_X)$. Then the spectral sequence of a double
complex (or an elementary substitute) yields our
conclusion.\end{proof} Morphisms $f$ satisfying the hypotheses of
the Theorem occur in diverse situations, e.g. regular fibre spaces
and resolutions of rational singularities. The Theorem says that
those schemes $X$ which admit a structure such as $f$ form an open
subset of of the moduli of $X$.
\begin{thm} Let $f:X\to Y$ be a proper surjective morphism \'etale
in codimension 1 where $X$ is normal. Then any deformation of $f$ is
determined by the associated deformation of $Y$.
\end{thm}\begin{proof} It will suffice to prove that the mapping
cone $\t_X\to f^!\t_Y$ is exact in degrees $\leq 1$. Working
locally, we may view $X$ as a subscheme of $Y\times R$, $R$ an
affine space, and consider a free resolution $J^.$ of the ideal of
$X$ in $Y\times R$. Let $K$ be the kernel of the natural map
$\hom(J^0,\O_X)\to\hom(J^{-1}, \O_X).$ Then $K$ is a torsion-free
$\O_X$ module and the natural map $k:T_R\otimes\O_X\to K$ is an
isomorphism in codimension 1 by our assumption that $f$ is \'etale
in codimension 1. As $X$ is normal, this easily implies $k$ is an
isomorphism, which proves our assertion.
\end{proof} It follows, e.g. that a small resolution of a
singularity $Y$ is locally uniquely determined by $Y$, i.e. cannot
be deformed without deforming $Y$. Smallness is of course essential
here. Under stronger hypotheses on the size of the exceptional
locus, we can actually identify the deformations of $f$ and $Y$
(compare \cite{Rhol}, Thm. 3.5):\begin{thm} Let $f:X\to Y$ be a
proper surjective morphism \'etale in codimension 2 with $X$ smooth.
Then any deformation of $Y$ lifts to a deformation of
$f$.\end{thm}\begin{proof} The idea is to relate deformations to
K\"ahler differentials, e.g. relate $H^i(\t_Y)$ to $\Ext^i(\Omega_Y,
\O_Y)$, so that we may apply Ischebeck's Lemma (
\cite{isch},\cite{mat} p. 104).
We work in the affine setting, with $Y\to P$ an affine embedding.
Let $F^.$ be a free resolution of $I=I_{Y/P}$ and consider the
complex
\begin{equation}\hat{\Omega}_Y:...\to F^{-1}\otimes\O_Y\to
F^0\otimes\O_Y\to\Omega_P\otimes\O_Y\end{equation} where the last
map corresponds to the derivative map $I\to\Omega_P\otimes\O_Y.$ By
definition, applying $\hom(.,\O_Y)$ to this complex yields the
tangent \sela $\t_Y.$ Then $\hat{\Omega}_Y$ induces a complex of
free $\O_X$-modules
\begin{equation}f^!\hat{\Omega}_Y:...\to F^{-1}\otimes\O_X\to
F^0\otimes\O_X\to\Omega_P\otimes\O_X\end{equation} and applying
$\Hom^.(\_\ ,\O_X)$ to this yields $f^!\t_Y.$ Now let $J^.$ be an
$\O_Y$-free resolution of $I/I^2,$ and consider the
complex\begin{equation}\tilde{\Omega}_Y:...\to J^{-1}\to
J^0\to\Omega_P\otimes\O_Y\end{equation}  We may assume that
\begin{equation}J^i\simeq F^i\otimes\O_Y, i=0,1;\ \
F^{-2}\otimes\O_Y\hookrightarrow J^{-2}.\end{equation} We have a
natural map $\hat{\Omega}_Y\to\tilde{\Omega}_Y$ hence
$f^!\hat{\Omega}_Y\to f^!\tilde{\Omega}_Y$, and these induce maps
\begin{equation}\label{t-omega} \Ext^i(\tilde{\Omega}_Y, \O_Y)\to
H^i( \t_Y), \Ext^i(f^!\tilde{\Omega}_Y, \O_X)\to
H^i(f^!\t_Y)\end{equation} that are bijective for $i\leq 1$,
injective for $i=2$. Note that
$H^0(f^!\tilde{\Omega}_Y)=f^*\Omega_Y$, while for $i<0,
H^i(f^!\tilde{\Omega}_Y)$ is supported on the exceptional locus of
$f$. By Ischebeck, it follows that
\begin{equation}\label{omega-omega} \Ext^i(f^!\tilde{\Omega}_Y,
\O_X)=\Ext^i(f^*\Omega_Y, \O_X), i\leq 2.\end{equation} By Ischebeck
again, the natural map
\begin{equation}\label{x-y}\Ext^i(\Omega_X, \O_X)=H^i(\t_X)\to
\Ext^i(f^*\Omega_Y, \O_X)\end{equation} is bijective for $i\leq 1$
and injective for $i=2$, because the cokernel of
$f^*\Omega_Y\to\Omega_X$ is supported on the exceptional locus of
$f$.\par Now given (\ref{t-omega}-\ref{x-y}), the completion of the
proof is boilerplate. Thus suppose given a 1st order deformation
$\alpha$ of $Y$. By (\ref{t-omega}) $f^!\alpha$ yields an element of
$\Ext^1(f^!\tilde{\Omega}_Y, \O_X)$ hence by (\ref{omega-omega}) an
element of $\Ext^1(f^*\Omega_Y, \O_X),$ which by (\ref{x-y}) comes
from an element $\beta$ of $H^1(\t_X)$, i.e. a 1st order deformation
of $X$ compatible with $\alpha$. Thus, $\alpha$ lifts to a
deformation of $f$. If $\alpha$ lifts to 2nd order, a suitable
obstruction element vanishes, and applying a similar argument to the
obstruction of $\beta$ then shows that this obstruction must vanish
too. Continuing in this manner, we show that any deformation of $Y$,
of any order, lifts to a deformation of $f$.
\end{proof} The Theorem does not extend to morphisms with a
codimension-2 exceptional locus, such as a small resolution of a
3-fold ODP.

\end{document}